\documentclass[11pt, a4paper]{article}

\usepackage{amsmath, amssymb, amsthm}
\usepackage{color}

\newcommand{\R}{\mathbb{R}}

\newcommand{\nor}{\parallel}

\newcommand{\bq}{\begin{equation}}
\newcommand{\bqw}{\begin{equation*}}
\newcommand{\eq}{\end{equation}}
\newcommand{\eqw}{\end{equation*}}
\newcommand{\vo}{\mathrm{dvol}}
\newcommand{\ov}{\overline}
\newcommand{\lm}{\lambda^+_{min}}
\newcommand{\vol}{\mathrm{vol}}
\newcommand{\supp}{\mathrm{supp\ }}
\newcommand{\grad}{\mathrm{grad}}
\renewcommand{\i}{\mathrm{i}}

\renewcommand{\Re}{\mathrm{Re}}
\newcommand{\Spin}{\mathrm{Spin}}
\newcommand{\SO}{\mathrm{SO}}

\newtheorem{lemma}{Lemma}[section]
\newtheorem{corollary}[lemma]{Corollary}
\newtheorem{theorem}[lemma]{Theorem}
\theoremstyle{definition}
\newtheorem{example}[lemma]{Example}
\newtheorem{definition}[lemma]{Definition}
\newtheorem{remark}[lemma]{Remark}

\title{The Hijazi inequality on conformally parabolic manifolds}
\author{Nadine Gro\ss e }
\date{}
\begin{document}

\maketitle
\begin{abstract}We prove the Hijazi inequality, an estimate for Dirac eigenvalues, for complete manifolds of finite volume. Under some additional assumptions on the dimension and the scalar curvature, this inequality is also valid for elements of the essential spectrum. This allows to prove the conformal version of the Hijazi inequality on conformally parabolic manifolds if the spin analog to the Yamabe invariant is positive.
\end{abstract}
\footnotetext{ 
keywords: Dirac operator -- conformally parabolic manifold -- conformal geometry
}

%\author{Nadine Gro\ss e %\footnote {I am thankful to Hans-Bert Rademacher for discussions and help.}
%}

%%%%%%%%%%%%%%%%%%%%%%%%%%%%%%%%%%%%%
\section{Introduction}
%%%%%%%%%%%%%%%%%%%%%%%%%%%%%%%%%%%%%

On a closed $n$-dimensional Riemannian spin manifold $(M,g,\sigma)$ with scalar curvature $s_g$, Friedrich \cite[Thm. A]{Frie80} gave an estimate for an eigenvalue $\lambda$ of the classical Dirac operator $D_g$:
\[ \lambda^2\geq \frac{n}{4(n-1)}\inf_M s_g. \]
This inequality was improved by Hijazi \cite{Hija} for dimension $n\geq 3$
\[\lambda^2\geq \frac{n}{4(n-1)}\mu,\]
where $\mu$ is the smallest eigenvalue of  $L_g\!=\!4\frac{n-1}{n-2}\Delta_g+s_g$, the conformal Laplacian.

On closed manifolds, there is a conformal version of the Hijazi inequality that relates the corresponding conformal quantities, that means the Yamabe invariant

 \[Q(M,g)=\inf \Bigg\{ \int\limits\limits_M vL_gv\vo_g\ \Big|\ \Vert v\Vert_{\frac{2n}{n-2}}=1, v\in C_c^\infty (M)\Bigg\}\]

with the $\lm$-invariant

\[ \lm (M,g,\sigma) =\inf_{g_0\in [g],\ \vol(M,g_0)<\infty} \lambda_1^+(M,g_0,\sigma)\vol (M,g_0)^{\frac{1}{n}}\]

where 
\bqw \lambda_1^+(M, 
g, \sigma)\!=\!\inf\Bigg\{\frac{\Vert D_g\phi\Vert^2}{(D_g\phi,\phi)}\ \Bigg
\vert\ 0\!<\!(D_g\phi,\phi )%\!:=\!\!\int\limits_M\!\!\!<\!D\phi,\! \phi\!>\!\vo_g
,\ \phi\in C_c^\infty(M,S)\Bigg\}\eqw

and $[g]$ is the set of all metrics conformal to $g$.  Furthermore, $\Vert .\Vert:=\Vert .\Vert_{L^2}$ and $C_c^\infty(M,S)$ denotes the compactly supported smooth spinors on $(M,g,\sigma)$ (Clearly, on closed manifolds $C_c^\infty(M,S)=C^\infty (M,S)$ is the set of all smooth spinors. But in order to use this definition later on noncompact manifolds as well we already wrote $C_c^\infty(M,S)$ here.).

The conformal Hijazi inequality
reads \[ \lm(M,g,\sigma)^2\geq \frac{n}{4(n-1)}Q(M,g).\]
This can be seen immediately
since on closed manifolds $\lambda_1^+$ is just the lowest positive Dirac eigenvalue, and for $Q\geq 0$ we have 
\bq\label{Qinf} Q(M,g)=\inf_{g_0\in[g],\  \vol(M,g_0)<\infty} \mu(g_0)\vol (M,g_0)^{\frac{2}{n}}\eq
where $\mu(g_0)$ is the infimum of the spectrum of $L_{g_0}$.

We note, that the $\lm$-invariant can also be defined as a variational problem similar to the Yamabe invariant \cite{Ammha}.\\
 
Since both the Yamabe and the $\lm$-invariant can also be considered on open manifolds, cf. \cite{Schoe}, \cite{NG}, it is interesting to know whether the conformal Hijazi inequality also holds on these manifolds.

In this paper, we examine this question for conformally parabolic manifolds, i.e. those that admit a complete metric of finite volume in their conformal class.

At first, we obtain an Hijazi equality for Riemannian manifolds equipped with a complete metric of finite volume. 

\begin{theorem}\label{Hij1}
Let $(M,g,\sigma)$ be a complete Riemannian spin manifold of finite volume and dimension $n>2$. Moreover, let $\lambda$ be an eigenvalue of its Dirac operator $D_g$, and let $\mu$ be the infimum of the spectrum of the conformal Laplacian. 
Then the following inequality holds:
\[ \lambda^2\geq\frac{n}{4(n-1)}\mu. \]
If equality is attained for a nonzero $\lambda$, the manifold admits a real Killing spinor and, hence, has to be Einstein and closed.
\end{theorem}

On complete manifolds, the Dirac operator is essentially self-adjoint and, in general, its spectrum consists of eigenvalues and the essential spectrum. (Note that also the Laplacian and, thus, the conformal Laplacian are formally self-adjoint). For elements of the essential spectrum, we also obtain an Hijazi-type inequality:

\begin{theorem}\label{n5sb}
Let $(M,g,\sigma)$ be a complete Riemannian spin manifold of dimension $n\geq 5$ with finite volume. Furthermore, let the scalar curvature of $(M,g)$ be bounded from below and $Q\ne 0$. If $\lambda$ is in the essential spectrum of the Dirac operator $\sigma_{ess}(D_g)$, then \[\lambda^2\geq \frac{n}{4(n-1)}\mu\]
where $\mu$ is the infimum of the spectrum of the conformal Laplacian
\end{theorem}

These two Hijazi inequalities allow to prove the conformal version:

\begin{theorem}\label{Hij2}
Let $(M,g, \sigma)$ be a conformally parabolic Riemannian spin manifold of dimension $n>2$. Let one of the following assumptions be fulfilled:
\begin{itemize}
\item[(i)] There is a complete conformal metric $\ov{g}$ of finite volume such that $0\not\in \sigma_{ess} (D_{\ov{g}})$.
\item[(ii)] The dimension is $n\geq 5$, and there is a complete conformal metric $\ov{g}$ of finite volume whose scalar curvature is bounded from below.
\end{itemize}
Then the conformal Hijazi inequality holds:
\[ \lm (M,g,\sigma)^2\geq \frac{n}{4(n-1)}Q(M,g).\]
\end{theorem}

A manifold that does not fulfill the assumption (i) has a vanishing $\lm$-invariant, cf. Lem. \ref{zero}.iii. Thus, we obtain

\begin{corollary}
Let $(M,g, \sigma)$ be a conformally parabolic Riemannian spin manifold of dimension $n>2$ and with $\lm>0$. Then the conformal Hijazi inequality is valid.
\end{corollary}

We give a brief outline of the paper:
In Section \ref{2}, we review some notations for the identification of spinor bundles of conformally equivalent metrics. Furthermore, we give a refined Kato inequality that will be used in the proof of Theorem \ref{n5sb}.
In Section \ref{3}, we give some properties of $\lambda_1^+$ and $\lm$ on conformally parabolic manifolds. With these preparations, the Theorems \ref{Hij1}, \ref{n5sb} and \ref{Hij2} can be shown in Section \ref{4}.\\

After this paper was written an article by Christian B\"ar \cite{Baer09} was published also dealing with the spectrum of noncompact manifolds. He restricts to uniformly positive curvature endomorphisms but also deals with generalized Dirac operator.
 
%%%%%%%%%%%%%%%%%%%%%%%%%%%%%%%%%%%%%
\section{Preliminaries}\label{2}
%%%%%%%%%%%%%%%%%%%%%%%%%%%%%%%%%%%%%

In this section, we first review the identification of spinor bundles of conformally equivalent metrics to fix notations. Then we give the refined Kato inequality that we use to prove Theorem \ref{n5sb}.\\

{\bf Spinor bundles of conformally related metrics}\nopagebreak

Let $\ov{g}=f^2g$ with $0<f\in C^\infty (M)$. Having fixed a spin structure $\sigma$ on $(M,g)$ with corresponding spinor bundle $S_g$, there always exists a corresponding spinor bundle $S_{\ov{g}}$ on $(M,\ov{g})$  and a vector bundle isomorphism \[A:\ S_g\to S_{\ov{g}},\ \psi\mapsto \ov{\psi}:=A(\psi)\] that is fibrewise an isometry \cite[Sect. 4.1]{Hija}, i.e. $\langle \ov{\psi}, \ov{\psi}\rangle_{\ov{g}}=\langle \psi,\psi\rangle_g$ or $|\ov{\psi}|_{\ov{g}}=|\psi|_g$. In the following, we write for both spinor bundles just $S$ and both norms just $|.|$ (which of the norms is meant follows from the inserted spinor).

Using this isometry, it is possible to compare the corresponding Dirac operators $D:=D_{g}$ and $\ov{D}:=D_{\ov{g}}$ \cite[Prop. 4.3.1]{Hija}:
\bqw \label{Dtrans}\ov{D}( f^{-\frac{n-1}{2}}\ov{\psi})=f^{-\frac{n+1}{2}}\ov{D\psi}.\eqw

{\bf Refined Kato inequalities}

The Kato inequality states that for any section $\phi$ of a Riemannian or Hermitian vector bundle $E$ endowed with a metric connection $\nabla$ on a Riemannian manifold $(M,g)$ we have 
\bq  2|\phi| |d|\phi||=|d|\phi|^2|=2| \langle \nabla\phi,\phi\rangle|\leq 2|\phi||\nabla\phi|, \label{Kato}\eq

i.e. $|d|\phi||\leq |\nabla\phi|$ away from the zero set of $\phi$. 

In \cite{CGH}, refined Kato inequalities were obtained for sections in the kernel of first-order elliptic differential operators $P$. They are of the form 
\[ |d|\phi||\leq k_P|\nabla\phi|\]
where $k_P$ is a constant $\leq 1$ depending on the operator $P$.\\

We sketch the set-up used in \cite{CGH}: Let $E$ be an irreducible natural vector bundle $E$ over an $n$-dimensional Riemannian (spin) manifold $(M,g)$ with fibrewise product $\langle .,.\rangle$ and a metric connection $\nabla$. Irreducible natural means that the vector bundle is obtained either from the %direct 
orthonormal frame bundle of M or from the spinor frame bundle with an irreducible representation of $\SO(n)$ or $\Spin(n)$ on a vector space $V$. We will denote this representation by $\lambda$. Further, let $\tau$ be the standard representation of $\SO(n)$ or $\Spin(n)$ on $\R^n$. Then the real tensor product $\tau\otimes\lambda$ splits into irreducible components as 

\[ \tau\otimes\lambda=\bigoplus_{j=1}^N\mu^j,\quad\quad\quad \R^n\otimes V=\bigoplus_{j=1}^NW_j.\]

This induces a decomposition of $T^*M\otimes E$ into irreducible subbundles $F_j$ associated to $\mu^j$. Further, let $\Pi_j$ denote the projection onto the $j$th summand of $\R^n\otimes V$ and $T^*M\otimes E$, respectively.

Let $P$ be a first-order linear differential operator of the form $P=\sum_{i\in I} \Pi_i\circ\nabla$ where $I\subseteq \{1,\ldots, N\}$. Moreover, we denote $\Pi_I:=\sum_{i\in I} \Pi_i$ and $\hat{I}:=\{1,\ldots,N\}\setminus I$.

Following the ansatz for the refined Kato inequalities we obtain the estimate:

\begin{lemma}\label{reKI2}{\rm \cite[Lem. 1.7.1]{NGDiss}} Let $P$ be an operator as defined above. Then we have away from the zero set of $\phi$
\[ |d|\phi||\leq |P\phi|+k_P|\nabla\phi|\]
where $k_P:=\sup_{|\alpha|=|v|=1} |\Pi_{\hat{I}} (\alpha\otimes  v)|$.
\end{lemma}

The proof can done analogously to the one of \cite{CGH} without the assumption that $\phi\in \ker P$. That's why the additional summand $|P\phi|$ appears and why the constant $k_P$ remains the same.

\begin{proof}
Let $\phi$ be a section of $E$. Then away from the zero set of $\phi$ we obtain

\begin{align*}
|d|\phi||&=\frac{|d|\phi|^2|}{2|\phi|}=\frac{|<\nabla\phi,\phi>|}{|\phi|}
\end{align*}
Let now $\alpha_0$ be a unit $1$-form with $<\nabla\phi,\phi>=c\alpha_0$ for some $c\in\R$. Then we have
\begin{align*}
<\nabla\phi,\alpha_0\otimes\phi>&=\sum_i <\nabla_{e_i}\phi,\alpha_0(e_i)\phi>=\sum_i\frac{1}{c}<\nabla_{e_i}\phi,\phi>^2\\
&=\sum_i\frac{<\nabla_{e_i}\phi,\phi>^2}{|<\nabla\phi,\phi>|}=|<\nabla\phi,\phi>|.
\end{align*}
Thus, we obtain
\begin{align*}
|d|\phi||&=\frac{|<\nabla\phi,\alpha_0\otimes\phi>|}{|\phi|}\\
&=\frac{|<(\Pi_I+\Pi_{\hat{I}})\nabla\phi,\alpha_0\otimes\phi>|}{|\phi|}\\
&\leq \frac{|<P\phi,\alpha_0\otimes\phi>|}{|\phi|}+\frac{|<\nabla\phi,\Pi_{\hat{I}}(\alpha_0\otimes\phi)>|}{|\phi|}\\
&\leq |P\phi|+|\nabla\phi|\sup_{|\alpha|=|v|=1} |\Pi_{\hat{I}}(\alpha\otimes v)|=|P\phi|+k_P|\nabla\phi|.
\end{align*}
\end{proof}

For the shifted (classical) Dirac operator $D-\lambda$ we have $k=\sqrt{\frac{n-1}{n}}$ \cite[(3.9)]{CGH}.

%%%%%%%%%%%%%%%%%%%%%%%%%%%%%%%%%%%%%
\section{The invariant on conformally parabolic manifolds}\label{3}
%%%%%%%%%%%%%%%%%%%%%%%%%%%%%%%%%%%%%

Firstly, we give a characterisation of conformally parabolic manifolds and consider the example of the Euclidean space.\\
In the rest of this section, we provide some properties of $\lambda_1^+$ for complete metrics with finite volume. 

\begin{definition}\cite[Sect. 3]{ZK}\label{ZK}
A Riemannian manifold is conformally parabolic if and only if its conformal class contains a complete metric of finite volume.
\end{definition}

\begin{example}\label{rnconf} Let $(M^m, g_M)$ be a closed $m$-dimensional Riemannian manifold. Then $(M\times (1,\infty), g=g_M+dt^2)$ is conformally parabolic since the conformal metric $\overline{g}=\frac{1}{t^2}g$ is complete and of finite volume.\\
Furthermore, for the new metric and for dimension $n>2$ the scalar curvature is calculated as (where $h=t^{-\frac{n-2}{2}}$ and $n=m+1$)

\begin{align*} s_{\ov{g}}&=4\frac{n-1}{n-2}h^{-\frac{n+2}{n-2}}\Delta h
+s_gh^{-\frac{4}{n-2}}\\
&=-4\frac{n-1}{n-2}t^{\frac{n+2}{2}}\left(1-\frac{n}{2}\right)\left(-\frac{n}{2}\right)t^{-\frac{n+2}{2}}+s_Mt^2\\
&= -(n-1)n+s_Mt^2.
\end{align*}

\end{example}

Next we give some properties of $\lambda_1^+$:

\begin{lemma}\label{zero}\hfill\\
\textbf{i)}\hphantom{ii}\  If $\lambda_1^+ (M,g, \sigma)=0$ and $\vol(M,g)<\infty$, then $\lm(M,g, \sigma)=0$.\\
\textbf{ii)}\ If $(M,g)$ is complete and $\lambda>0$ is an eigenvalue of $D$ or an element of its essential spectrum, then $\lambda_1^+(M,g,\sigma)\leq\lambda$.\\
\textbf{iii)}\ A complete Riemannian spin manifold of finite volume for which there exists $\lambda>0$ in the essential spectrum of its Dirac operator has a vanishing $\lm$-invariant.
\end{lemma}

\begin{proof} i) is seen immediately from the definition of $\lm$.\\
ii) There exists a sequence $\phi_i\in C_c^\infty(M,S)$ with $\Vert D\phi_i-\lambda\phi_i\Vert\to 0$ and $\Vert\phi_i\Vert\to 1$: If $\lambda$ is in the essential spectrum, this follows directly from the definition. If $\lambda$ is an eigenvalue with eigenspinor $\phi\in C^\infty(M,S)\cap L^2(M,S)$, we choose $\phi_i=\eta_i\phi$ where $\eta_i$ is a smooth cut-off function such that $\eta_i\equiv 1$ on $B_i(p)$ ($p\in M$ fixed), $\eta_i\equiv 0$ on $M\setminus B_{2i}(p)$ and in between $|\nabla \eta_i|\leq\frac{2}{i}$. This is always possible since $(M,g)$ is complete. Then $\phi_i$ is the sequence in demand since $\Vert (D-\lambda)\phi_i\Vert=\Vert\nabla\eta_i\cdot\phi\Vert\leq \frac{2}{i}\Vert\phi\Vert$.\\
Thus, in both cases \[\frac{\Vert D\phi_i\Vert^2}{(D\phi_i,\phi_i)}\to \lambda\]
which proves the claim.\\
iii) Since the essential spectrum only depends on the manifold at infinity, see \cite[Prop. 1]{Ba3}, there is a sequence $\phi_i\in C_c^\infty(M\setminus B_i(p),S)$ ($p\in M$ fixed) with $\Vert (D-\lambda)\phi_i\Vert\to 0$ and $\Vert\phi_i\Vert=1$. Thus, as in ii) we find 
\[\lm (M\setminus B_r(p),g, \sigma)\leq \lambda\, \vol(M\setminus B_r(p) ,g)^\frac{1}{n}\to 0\]
for $r\to \infty$. With $\lm(M,g,\sigma)\leq \lm (M\setminus B_r(p),g, \sigma)$ (\cite[Lem. 2.1]{NG} or it can be seen directly from the variational description of $\lm$ \cite{Ammha} using that every spinor compactly supported on $M\setminus B_r(p)$ is also compactly supported on $M$), we have $\lm (M,g,\sigma)=0$.
\end{proof}

\begin{lemma}\label{spec} Let $(M,g, \sigma)$ be a complete Riemannian spin manifold. Then 
\[\lambda_1^+(M,g,\sigma)=\inf \{\sigma(D)\cap (0,\infty)\}\]
where $\sigma(D)$ denotes the Dirac spectrum.\end{lemma}

\begin{proof} Since $(M,g)$ is complete, $D$ is essentially self-adjoint and has no residual spectrum, cf. \cite[Chapt. 4]{Frie}.
By the spectral theorem for unbounded self-adjoint operators, we obtain that for every $\phi\in C_c^{\infty}(M,S)$ with $(D\phi, \phi)>0$ 
\begin{align*} \frac{\Vert D\phi\Vert^2}{(D\phi,\phi)}&=\frac{\int_{\sigma (D)} \lambda^2\ d\langle E_\lambda \phi,\phi\rangle}{\int_{\sigma (D)} \lambda\ d\langle E_\lambda \phi,\phi\rangle }\geq \frac{\int_{\sigma (D)\cap (0,\infty)} \lambda^2\ d\langle E_\lambda \phi,\phi\rangle }{\int_{\sigma (D)\cap (0,\infty)} \lambda\ d\langle E_\lambda \phi,\phi\rangle }\displaybreak[0]\\
&\geq \frac{\lambda_0\int_{\sigma (D)\cap (0,\infty)} \lambda\ d\langle E_\lambda \phi,\phi\rangle }{\int_{\sigma (D)\cap (0,\infty)} \lambda\ d\langle E_\lambda \phi,\phi\rangle }=\lambda_0
  \end{align*}
where $\lambda_0=\inf \{\sigma(D)\cap (0,\infty)\}$. Note that the denominator \[\int_{\sigma (D)\cap (0,\infty)} \lambda\ d\langle E_\lambda \phi,\phi\rangle \geq \int_{\sigma (D)} \lambda\ d\langle E_\lambda \phi,\phi\rangle=(D\phi,\phi)>0.\] Hence, we have $\lambda_1^+\geq \inf \{\sigma(D)\cap (0,\infty)\}$.\\
The converse inequality is obtained by Lemma \ref{zero}.ii.\end{proof}

From Lemma \ref{zero}.iii and Lemma \ref{spec}, we have	

\begin{corollary}\label{poslam}
Let $(M,g,\sigma)$ be a complete Riemannian spin manifold of finite volume with $\lm >0$. Then $\sigma (D)\cap (0,\infty)$ consists only of eigenvalues. 
\end{corollary}

The next Lemma shows that for defining the $\lm$-invariant on conformally parabolic manifolds we do not need the infimum over all conformal metrics.

\begin{lemma}\label{sequ}
Let $(M,g,\sigma)$ be a conformally parabolic Riemannian spin manifold. Then there exists a sequence of complete conformal metrics $g_i$ of unit volume such that $\lambda_1^+(g_i)\to \lm (g)$ and $g_i\equiv g_1$ near infinity, i.e. 
\[ \lm (M,g,\sigma)=\inf\{ \lambda_1^+(M,\ov{g},\sigma)\ |\ \ov{g}\equiv g_1 \mathrm{\ near\ infinity},\ \vol (M,\ov{g})=1\},\] 
where ``near infinity'' refers to the existence of compact subsets $U_{\ov{g}}\subset M$ such that $\ov{g}\equiv g_1$ on $M\setminus U_{\ov{g}}$.
\end{lemma}

\begin{proof} Assume that $g=g_1$ is already complete and of unit volume. Let $g_i=f_i^2g$ be a sequence of conformal metrics of unit volume with $\lambda_1^+(g_i)\to \lm$ for $i\to \infty$. Thus, there is a sequence $\phi_i\in C_c^{\infty}(M,S)$ such that \[F(\phi_i,g_i):=\frac{\nor D_{g_i}\phi_i\nor_{g_i}^2}{(D_{g_i}\phi_i, \phi_i)_{g_i}}\to \lm\]
Now, we choose the conformal factor $h_i$ such that $h_i$ is equal to $f_i$ on the support of $\phi_i$, $h_i=1$ near infinity and $\int_M h_i^n\vo_g=1$. Then, $F(\phi_i, h_i^2g)=F(\phi_i, g_i)\to \lm$, the metrics $h_i^2g$ are complete since $g$ is complete, and they have unit volume.
\end{proof}

%%%%%%%%%%%%%%%%%%%%%%%%%%%%%%%%%%%%%
\section{Proof of Hijazi inequalities}\label{4}
%%%%%%%%%%%%%%%%%%%%%%%%%%%%%%%%%%%%%

Firstly, we follow the main idea of the proof of the original Hijazi inequality, but we fix the used conformal factor with the help of an eigenspinor. This results in a conformal metric on the manifold without the zero-set of the eigenspinor and we have to use cut-off functions near this zero-set and near infinity to obtain compactly supported test functions.
\begin{proof}[Proof of Theorem \ref{Hij1}]
Let $\psi\in C^\infty(M,S)\cap L^2(M,S)$ be an eigenspinor satisfying $D\psi=\lambda\psi$ and $\Vert \psi\Vert=1$. Its zero-set $\Omega$ is closed and contained in a closed countable union of smooth $(n-2)$-dimensional submanifolds which has locally finite $(n-2)$-dimensional Hausdorff measure \cite[p. 189]{Ba99}.

We fix a point $p\in M$. Since $M$ is complete, there exists a cut-off function $\eta_i: M\to [0,1]$ which is zero on $M\setminus B_{2i}(p)$ and one on $B_{i}(p)$. In between the function is chosen such that $|\nabla \eta_i|\leq \frac{4}{i}$ and $\eta_i\in C_c^\infty(M)$.

While $\eta_i$ cuts off $\psi$ at infinity, we define another cut-off near the zeros of $\psi$. For this purpose, we can assume without loss of generality that $\Omega$ is  itself the countable union of $(n-2)$-submanifolds with locally finite $(n-2)$-dimensional Hausdorff measure described above.

Let now $\rho_{a,\epsilon}$ be defined as
\[\rho_{a,\epsilon}(x)=\Bigg\{ \begin{array}{ll}  0  &\mathrm{for\ } r< a\epsilon\\ 
1-\delta\ln\frac{\epsilon}{r} &\mathrm{for\ } a \epsilon \leq r\leq\epsilon\\
1  &\mathrm{for\ } \epsilon < r \end{array}
 \]
where $r=d(x,\Omega)$ is the distance from $x$ to $\Omega$. The constant $a<1$ is chosen such that $\rho_{a,\epsilon} (a\epsilon)=0$, i.e. $a=e^{-\frac{1}{\delta}}$. For $\epsilon$ small enough $\rho_{a,\epsilon}$ is well-defined and $\rho_{a,\epsilon}$ is continuous and Lipschitz. Now, we define $\psi_{ia}:=\eta_i\rho_{a,\epsilon}\psi$. 
Sincer $\psi\in C^\infty (M,S)$ the spinor $\psi_{ia}$ is an element in $ H_1^r(M,S)$ for all $1\leq r\leq \infty$.

 These spinors are compactly supported on $M\setminus \Omega$. Furthermore, $\ov{g}=e^{2u}g=h^{\frac{4}{n-2}}g$ with $h=|\psi|^{\frac{n-2}{n-1}}$ is a metric on $M\setminus{\Omega}$. Setting $\ov{\phi_{ia}}:=e^{-\frac{n-1}{2}u}\ov{\psi_{ia}}$\ ($\phi= e^{-\frac{n-1}{2}u}\psi$), the Lichnerowicz-type formula \cite[(5.4)]{Hija} implies
\begin{align*}\Vert (\ov{D}-&\lambda e^{-u})\ov{\phi_{ia}}\Vert^2_{\ov{g}}=\Vert \ov{\nabla}^{\lambda e^{-u}}\ov{\phi_{ia}} \Vert^2_{\ov{g}}+ \int\limits_{M\setminus\Omega} \left(\frac{\ov{s}}{4}-\frac{n-1}{n}\lambda^2e^{-2u}\right)|\ov{\phi_{ia}}|^2\vo_{\ov{g}}\\
&\hphantom{=} -\frac{n-1}{n}(2\lambda e^{-u}(\ov{D}-\lambda e^{-u})\ov{\phi_{ia}}+\lambda e^{-u}\ov{\grad\, e^{-u}\cdot \phi_{ia}},\ov{\phi_{ia}})_{\ov{g}}\displaybreak[0]\\
&=\Vert \ov{\nabla}^{\lambda e^{-u}}\ov{\phi_{ia}} \Vert^2_{\ov{g}}+ \int\limits_M \left(\frac{\ov{s}}{4}-\frac{n-1}{n}\lambda^2e^{-2u}\right)e^u|{\psi_{ia}}|^2\vo_{{g}}\\
&\hphantom{=} -2\frac{n-1}{n}\Re((D-\lambda)\psi_{ia},\lambda e^{-u}\psi_{ia})_g\displaybreak[0]\\
&=\Vert \ov{\nabla}^{\lambda e^{-u}}\ov{\phi_{ia}} \Vert^2_{\ov{g}}+\frac{1}{4} \int\limits_M h^{-1}Lh\,e^{-u}|\psi_{ia}|^2\vo_g\displaybreak[0]\\
&\hphantom{=}-\frac{n-1}{n}\lambda^2\int\limits_M e^{-u}|\psi_{ia}|^2 \vo_g -2\frac{n-1}{n}\Re((D-\lambda)\psi_{ia},\lambda e^{-u}\psi_{ia})_g,
\end{align*} 
where $\nabla^f_X\phi:=\nabla_X\phi+\frac{f}{n}X\cdot\phi$  for $f=\lambda e^{-u}\in C^\infty(M)$ is the Friedrich connection. For the second line we used $|\ov{\phi_{ia}}|^2\vo_{\ov{g}}=e^u|\psi_{ia}|^2\vo_g$ and that the term $(\lambda e^{-u}\ov{\grad e^{-u}\cdot \phi_{ia}},\ov{\phi_{ia}})_{\ov{g}}\in \i\R$ since $\langle \nabla f\cdot \phi, \phi\rangle \in \i\R$, cf. \cite[Lem. 3.1]{Hija}. The last line is obtained by replacing $\ov{s}e^{2u}=h^{-1}Lh$.

With $D\psi=\lambda\psi$ and $\langle \nabla f\cdot \psi, \psi\rangle \in \i\R$, we obtain \[\Re((D-\lambda)\psi_{ia},\lambda e^{-u}\psi_{ia})_g=\Re(\nabla (\eta_i\rho_{a,\epsilon})\psi, \lambda e^{-u}\eta_i\rho_{a,\epsilon}\psi)_g=0.\]  Inserting this result,  $\ov{D}\,\ov{\phi}=\lambda e^{-u}\ov{\phi}$ and $\Vert \ov{\nabla}^{\lambda e^{-u}}\ov{\phi_{ia}} \Vert^2_{\ov{g}}\geq 0$  into the formula from above and replace $e^u=|\psi|^{\frac{2}{n-1}}$ and $h=|\psi|^{\frac{n-2}{n-1}}$ we further have
\begin{align*}
\Vert \ov{\nabla} (\eta_i\rho_{a,\epsilon})\ov{\phi}\Vert^2_{\ov{g}}\geq\! \int\limits_M\!\left(\frac{1}{4}\eta_i^2\rho_{a,\epsilon}^2 |\psi|^{\frac{n-2}{n-1}} L|\psi|^{\frac{n-2}{n-1}}-\frac{n-1}{n}\lambda^2\eta_i^2\rho_{a,\epsilon}^2 |\psi|^{2\frac{n-2}{n-1}}\right)\!\vo_g.
\end{align*}

Moreover, we have
\begin{align*}\Vert \ov{\nabla} (\eta_i\rho_{a,\epsilon})\ov{\phi}\Vert^2_{\ov{g}}=\int\limits_M |e^{-u}\ov{\nabla (\eta_i\rho_{a,\epsilon})\cdot\phi}|^2\vo_{\ov{g}}=\int\limits_M |\nabla (\eta_i\rho_{a,\epsilon})\cdot\psi|^2e^{-u}\vo_{g}.
\end{align*}

Thus, with $e^u=|\psi|^\frac{2}{n-1}$ the above inequality reads
\begin{align*}
&\int\limits_M |\nabla (\eta_i\rho_{a,\epsilon})|^2 |\psi|^{2\frac{n-2}{n-1}}\vo_g  \geq \frac{1}{4}\int\limits_M \eta_i\rho_{a,\epsilon} |\psi|^{\frac{n-2}{n-1}} L (\eta_i\rho_{a,\epsilon}|\psi|^{\frac{n-2}{n-1}})\vo_g \\[0.3cm]
&\hphantom{aa}-\frac{n-1}{n-2}\int\limits_M |\nabla (\eta_i\rho_{a,\epsilon})|^2 |\psi|^{2\frac{n-2}{n-1}}\vo_g -\frac{n-1}{n}\lambda^2\int\limits_M\eta_i^2\rho_{a,\epsilon}^2 |\psi|^{2\frac{n-2}{n-1}}\vo_g.
\end{align*} 

Hence, we obtain
\[
\frac{2n-3}{n-2}\!\int\limits_M\!\!|\nabla (\eta_i\rho_{a,\epsilon})|^2 |\psi|^{2\frac{n-2}{n-1}}\vo_g  \geq \left(\frac{\mu}{4} -\frac{n-1}{n}\lambda^2\right)\!\int\limits_M\eta_i^2\rho_{a,\epsilon}^2 |\psi|^{2\frac{n-2}{n-1}}\vo_g,\]
where $\mu$ is the infimum of the spectrum of the conformal Laplacian.
With $(a+b)^2\leq 2a^2+2b^2$ we have 
\[ K\!\int\limits_M\!(\eta_i^2|\nabla \rho_{a,\epsilon}|^2 +\rho_{a,\epsilon}^2|\nabla \eta_i|^2) |\psi|^{2\frac{n-2}{n-1}}\vo_g  \geq \left(\frac{\mu}{4}-\frac{n-1}{n}\lambda^2\right)\!\Vert \eta_i\rho_{a,\epsilon} |\psi|^{\frac{n-2}{n-1}}\Vert_g^2\]
where $K=2\frac{2n-3}{n-2}$.

Now we let $a$ tend to zero (for fixed $i$):

Recall that $\Omega\cap \ov{B_{2i}(p)}$ is bounded, closed, $(n-2)$-$C^\infty$-rectifiable and has still locally finite $(n-2)$-dimensional Hausdorff measure. For fixed $i$ we estimate

\[ \int\limits_M |\nabla \rho_{a,\epsilon}|^2\eta_i^2|\psi|^{2\frac{n-2}{n-1}}\vo_g\leq \sup_{B_{2i}(p)} |\psi|^{2\frac{n-2}{n-1}}\ \int\limits_{B_{2i}(p)} |\nabla \rho_{a,\epsilon}|^2\vo_g.\]

For $y\in \Omega$ we set
$ B_\epsilon^{2}(y):=\{x\in B_\epsilon\ |\ d(x,y)=d(x,{\Omega})\}$ with $B_\epsilon:=\{x\in M\ |\ d(x,{\Omega})\leq\epsilon\}$. Since $\Omega$ is locally a $(n-2)$-dimensional submanifold, $B_\epsilon^2(y)$ is two-dimensional. Moreover, there is an inclusion $B_\epsilon^2(y)\hookrightarrow B_\epsilon(0)\subset\R^{2}$ via the normal exponential map. For $\epsilon$ sufficiently small the image of $ B_\epsilon^{2}(y)$ under this inclusion is star shaped for all $y\in \Omega\cap B_{2i}(p)$. Then we can calculate the 
\begin{align*}  \int\limits_{ B_{2i}(p)} |\nabla \rho_{a,\epsilon}|^2\vo_g&=\int\limits_{B_\epsilon\cap B_{2i}(p)}\hspace{-0.5cm} |\nabla \rho_{a,\epsilon}|^2\vo_g
= \int\limits_{x\in \Omega \cap B_{2i}(p)} \int \limits_{B^{2}_\epsilon(x)\setminus B^{2}_{a\epsilon}(x)}
\hspace{-0.5cm} |\nabla \rho_{a,\epsilon}|^2\vo_{g}\\
& \leq { \vol_{n-2}}
(\Omega\cap B_{2i}(p) ) \sup_{x\in\Omega\cap B_{2i}(p)} \int\limits_{B^{2}_\epsilon(x)\setminus B^{2}_{a\epsilon}(x)}
\hspace{-0.5cm} |\nabla \rho_{a,\epsilon}|^2\vo_{g_2}\\
&\leq c \vol_{n-2}(\Omega\cap B_{2i}(p) ) \int\limits_{B_\epsilon(0)\setminus {B}_{a\epsilon}(0)}\hspace{-0.5cm} |\nabla \rho_{a,\epsilon}|^2\vo_{g_E}\\
& \leq c^\prime  \int\limits^\epsilon_{a\epsilon} \frac{\delta^2}{r}dr=-c^\prime \delta^2\ln a=c^\prime\delta\to 0 \quad \mathrm{for\ } a\to 0
 \end{align*}
where $\vol_{n-2}$ denotes the $(n-2)$-dimensional volume %(Hausdorff measure?)
and $g_2=g_{|_{B_\epsilon^2(p)}}$. The positive constants $c$ and $ c^\prime$ arise from $\vol_{n-2}({\Omega}\cap B_{2i}(p))$ and the comparison of $\vo_{g_{2}}$ with the volume element of the Euclidean metric.

Furthermore, by the monotone convergence theorem, we obtain
\[\int\limits_{B_{2i}(p)}\hspace{-0.2cm} \rho_{a,\epsilon}^2 |\nabla \eta_i|^2 |\psi|^{2\frac{n-2}{n-1}}\vo_g\to \int\limits_{B_{2i}(p)}\hspace{-0.2cm} |\nabla \eta_i|^2 |\psi|^{2\frac{n-2}{n-1}}\vo_g\]
 as $a\to 0$ and, thus,
\[ K\int\limits_M |\nabla \eta_i|^2 |\psi|^{2\frac{n-2}{n-1}}\vo_g  \geq \left(\frac{\mu}{4}-\frac{n-1}{n}\lambda^2\right)\int\limits_M \eta_i^2 |\psi|^{2\frac{n-2}{n-1}}\vo_g.\\[0.2 cm] 
\]

Next, we want to establish the limit for $i\to \infty$:\\
Since $M$ has finite volume and $\Vert\psi\Vert=1$, the H\"older inequality ensures that $\int\limits_M |\psi |^{2\frac{n-2}{n-1}}\vo_g$ is bounded. With $|\nabla \eta_i|\leq \frac{4}{i}$ we obtain in the limit as $\i\to \infty$ 
\[ 0\geq \left(\frac{\mu}{4}- \frac{n-1}{n}\lambda^2\right)\int_M |\psi |^{2\frac{n-2}{n-1}}\vo_g\]
and, thus, 
\[ \lambda^2\geq \frac{n}{4(n-1)}\mu. \]

Equality is attained if and only if $\Vert \ov{\nabla}^{\lambda e^{-u}}\ov{\phi_{ia}} \Vert^2_{\ov{g}}\to 0$ for $i\to \infty$ and $a\to 0$. We have
\begin{align*} 0\leftarrow &\Vert \ov{\nabla}^{\lambda e^{-u}}\ov{\phi_{ia}} \Vert_{\ov{g}}=
\Vert \eta_i\rho_{a,\epsilon}\ov{\nabla}^{\lambda e^{-u}}\ov{\phi} +\ov{\nabla}(\eta_i\rho_{a,\epsilon}) \ov{\phi} \Vert_{\ov{g}}\\
&\geq \Vert \eta_i\rho_{a,\epsilon}\ov{\nabla}^{\lambda e^{-u}}\ov{\phi}\Vert_{\ov{g}} -{\Vert \ov{\nabla}(\eta_i\rho_{a,\epsilon}) \ov{\phi}\Vert_{\ov{g}}}.
\end{align*}
With $\Vert \ov{\nabla}(\eta_i\rho_{a,\epsilon}) \ov{\phi}\Vert_{\ov{g}}\to 0$, see above, 
$\ov{\nabla}^{\lambda e^{-u}}\ov{\phi}$ has to vanish on $M\setminus {\Omega}$. With \cite[Cor. 3.6]{Hija} this implies that $e^{-u}$ is constant. Thus, $(M,g)$ is Einstein and possesses a real Killing spinor, cf. \cite[p. 118]{Frie}. Moreover, if $\lambda>0$ its Einstein constant is positive. Thus, the Ricci curvature is a positive constant and, hence, due to the Theorem of Bonnet-Myers $M$ is already closed.
\end{proof}

Next, we want to prove Theorem \ref{n5sb} using the refined Kato inequality. Similar methods were used by Davaux in \cite{Dav03}. But before we state the following lemma:

\begin{lemma}\cite[Lem. 1.2.10 and 1.2.11]{NGDiss}\label{lem_ess}
 Let $(M,g,\sigma)$ be a complete Riemannian spin manifold. \\
 i) Let $\lambda \in\R$. Then $(D-\lambda)$ and $(D-\lambda)^2$ are essentially self-adjoint.\\
 ii) $0$ is in the essential spectrum of $D-\lambda$ if and only if $0$ is in the essential spectrum of $(D-\lambda)^2$.
If $0$ is in the essential spectrum of $D-\lambda$, then there is a normalized sequence $\phi_i\in C_c^\infty(M,S)$ such that $\phi_i$ converges $L^2$-weakly to $0$ and $\Vert (D-\lambda)\phi_i\Vert\to 0$ and $\Vert (D-\lambda)^2\phi_i\Vert\to 0$.
\end{lemma}

\begin{proof}
i) Since $D$ is essentially self-adjoint, $(D-\lambda)$ and $-2\lambda D$ are also essentially self-adjoint. From the inequality (see \cite[Prop. 6.2]{Wolf})
\bq\label{relbd} \Vert D\phi\Vert^2\leq t\Vert D^2\phi\Vert^2+\frac{1}{t}\Vert\phi\Vert^2 \quad\rm{for} t>0\eq we see that $D$ is $D^2$-bounded with relative bound $\sqrt{t}$ (For a definition of relative boundedness see \cite[Sect. X.2]{RS2}). Similarly, $-2\lambda D$ is $D^2$-bounded with relative bound $\frac{1}{2\lambda}\sqrt{t}$.\\
Then the Kato-Rellich Theorem \cite[Thm. X.12]{RS2} yields that $D^2-2\lambda D$ and, therefore, $(D-\lambda)^2$ is essentially self-adjoint.\\
ii) Due to i) both operators $A:=D-\lambda$ and $A^2=(D-\lambda)^2$ are essentially self-adjoint on $C_c^\infty(M,S)$.
Denote by $E_A$ and $E_{A^2}$ the projector-valued measures belonging to $A$ and $A^2$, respectively. We have $\supp E_{A^2}=[0,\infty)$ and   $E_{A^2}([a,b])=E_{A}([-\sqrt{b},-\sqrt{a}])+ E_{A}([\sqrt{a},\sqrt{b}])$ for $0\leq a\leq b$ which follows from \cite[Thm. 3.1]{BSU2}. Thus, if $0$ is in the (not necessarily essential) spectrum of $A$, then it is also contained in the spectrum of $A^2$ and vice versa.\\
Let now $0\in \sigma_{ess}(A)$. Then for every $\epsilon>0$ we obtain for the dimension of the image space of the projector $E_A([-\epsilon, \epsilon])$ that $\mathrm{dim} E_A([-\epsilon, \epsilon])H=\infty$ where $H:=L^2(M,S)$ and, thus, $\mathrm{dim} E_{A^2}([0,\epsilon^2])H=\infty$. Hence, we have $0\in \sigma_{ess}(A^2)$. Analogously, it follows from $0\in \sigma_{ess}(A^2)$, that $0\in \sigma_{ess}(A)$.\\
Next, let $0\in \sigma_{ess}((D-\lambda)^2)$. Due to the definition of $\sigma_{ess}$ there is a normalized sequence $\phi_i\in C_c^\infty(M,S)$ such that $\phi_i$ converges $L^2$-weakly to $0$ and $\Vert (D-\lambda)^2\phi_i\Vert\to 0$. Then, we have
\[ \Vert (D-\lambda)\phi_i\Vert^2=((D-\lambda)^2\phi_i,\phi_i)\leq \Vert(D-\lambda)^2\phi_i\Vert\,\Vert\phi_i\Vert\to 0.\]
\end{proof}

\begin{proof}[Theorem \ref{n5sb}] We may assume $\vol(M,g)=1$. If $Q>0$, there exists a compactly supported function $v$ with $\int_M vL_gv\vo_g <0$. Then $\mu<0$ as well and the claimed inequality is trivially fulfilled. Thus, we assume from now on that $Q>0$.

If $\lambda$ is in the essential spectrum of $D$, then $0$ is in the essential spectrum of $D-\lambda$. By Lemma \ref{lem_ess} there is a sequence $\phi_i\in C_c^\infty (M,S)$ such that $\Vert (D-\lambda)^2\phi_i\Vert\to 0$ and $\Vert(D-\lambda)\phi_i\Vert\to 0$ while $\Vert \phi_i\Vert=1$. We may assume that $|\phi_i|\in C_c^\infty(M)$. That can always be achieved by a small perturbation.

Now let $\frac{1}{2}\leq \beta\leq 1$. Then $|\phi_i|^\beta\in H_1^2(M)$.

Firstly, we will show that the sequence $\Vert d |\phi_i|^\beta\Vert$ is bounded:

By the H\"older inequality we have
\begin{align*}
0&\leftarrow \Vert\phi_i\Vert^{2\beta-1} \Vert (D-\lambda)^2\phi_i\Vert \geq \Vert |\phi_i|^{2\beta-1}\Vert_{\{|\phi_i|\ne 0\}} \Vert (D-\lambda)^2\phi_i\Vert \\
&\geq \Big|\int\limits_{|\phi_i|\ne 0} |\phi_i|^{2\beta-2}\langle(D-\lambda)^2\phi_i,\phi_i\rangle\vo_g\Big|.\end{align*}
Using the Lichnerowicz formula \cite[(5.4)]{Hija} where $\Delta^\lambda=(\nabla^{\lambda})^*\nabla^{\lambda}$, we obtain 
\begin{align*}
&\Vert (D-\lambda)^2\phi_i\Vert\\
&\geq \Bigg|\int\limits_{|\phi_i|\ne 0} |\phi_i|^{2\beta-2}\langle \Delta^\lambda\phi_i,\phi_i\rangle \vo_g %&\phantom{=}
+\int\limits \left(\frac{s}{4}-\frac{n-1}{n}\lambda^2\right)|\phi_i|^{2\beta}\vo_g
\\&\phantom{=}
-2\frac{n-1}{n}\int\limits_{|\phi_i|\ne 0} |\phi_i|^{2\beta-2} \langle (D-\lambda)\phi_i,\lambda\phi_i\rangle \vo_g\Bigg|\\
&\geq \int\limits_{|\phi_i|\ne 0}\hspace{-0.3cm} |\phi_i|^{2\beta-2}|\nabla^\lambda\phi_i|^2\vo_g+
%\!\\\phantom{=}+
2(\beta-1) \int\limits_{|\phi_i|\ne 0}\hspace{-0.3cm} |\phi_i|^{2\beta-3}\!\langle d|\phi_i|\!\cdot\!\phi_i,\nabla^\lambda\phi_i \rangle \vo_g\\
&\phantom{=}
+\!\int\!\left(\frac{s}{4}-\frac{n-1}{n}\lambda^2\right) |\phi_i|^{2\beta}\vo_g
-2\frac{n-1}{n}\lambda\Vert |\phi_i|^{2\beta-1}\Vert_{\{|\phi_i|\ne 0\}} \Vert(D-\lambda)\phi_i\Vert
\end{align*}
For the second summand on the right handside we then use $2(\beta-1)\leq 0$ and $\langle d|\phi_i| \cdot \phi_i,\nabla^\lambda\phi_i \rangle \leq |d|\phi_i|| \cdot |\phi_i| |\nabla^\lambda\phi_i|$ and for the first summand we use the ordinary Kato inequality $|\nabla ^{\lambda}\phi_i|\geq |d|\phi_i||$ (see \eqref{Kato}), to obtain
\begin{align*}
\Vert (D-&\lambda)^2\phi_i\Vert \geq (2\beta-1)\hspace{-0.2cm}\int\limits_{|\phi_i|\ne 0}\hspace{-0.2cm} |\phi_i|^{2\beta-2}|d|\phi_i|||\nabla^\lambda \phi_i|\vo_g\\
&+ \!\int\limits\left(\frac{s}{4}-\frac{n-1}{n}\lambda^2\right) |\phi_i|^{2\beta}\vo_g
\phantom{=}-2\frac{n-1}{n}\lambda\Vert \phi_i\Vert^{2\beta-1} \Vert(D-\lambda)\phi_i\Vert\end{align*}
Using $\beta \geq \frac{1}{2}$, again the Kato inequality and $\Vert \phi_i\Vert=1$ we get
\begin{align}
0&\leftarrow \Vert (D-\lambda)^2\phi_i\Vert\nonumber \\
&\geq (2\beta-1)\int\limits_{|\phi_i|\ne 0} |\phi_i|^{2\beta-2}|d|\phi_i||^2\vo_g+ \!\int\limits\left(\frac{s}{4}-\frac{n-1}{n}\lambda^2\right) |\phi_i|^{2\beta}\vo_g\nonumber\\
&\phantom{=}-2\frac{n-1}{n}\lambda \Vert(D-\lambda)\phi_i\Vert\nonumber\\
&\geq (2\beta-1)\frac{1}{\beta^2}\int\limits_{|\phi_i|\ne 0} |d|\phi_i|^\beta|^2\vo_g+ \int\limits\left(\frac{s}{4}-\frac{n-1}{n}\lambda^2\right) |\phi_i|^{2\beta}\vo_g\nonumber\\
&\phantom{=}-2\frac{n-1}{n}\lambda \Vert(D-\lambda)\phi_i\Vert\label{non-zero-in}
\end{align}

Since $s$ is bounded from below, $\int\limits s|\phi_i|^{2\beta}\vo_g\geq \inf s\ \Vert \phi_i \Vert^{2\beta}_{2\beta}\geq \min\{\inf s, 0\} $ is also bounded from below. Thus, with $\Vert (D-\lambda)\phi_i\Vert\to 0$ we see that $\Vert d|\phi_i|^\beta\Vert$ is also bounded.

With this preparation we can now prove the Hijazi inequality. We fix $\alpha=\frac{n-2}{n-1}$ and obtain 
\begin{align}
\Bigg(&\frac{\mu}{4}-\frac{n-1}{n}\lambda^2\Bigg)\Vert |\phi_i|^{\alpha}\Vert^2\nonumber\\
& \label{ineq-ess}\leq  \frac{1}{4}\int\limits |\phi_i|^\alpha L|\phi_i|^\alpha\vo_g-\frac{n-1}{n}\lambda^2\Vert|\phi_i|^\alpha\Vert^2\\
&=\!\int\!\! |\phi_i|^{2\frac{n-2}{n-1}-2}\left(\frac{n}{n-1} |d|\phi_i||^2+\frac{1}{2}d^*d|\phi_i|^2\nonumber
%\\&\phantom{aa}
+\left(\frac{s}{4}-\frac{n-1}{n}\lambda^2\right)\!|\phi_i|^2\right)\vo_g
\end{align}
where we used the definition of $\mu$ as infimum of the spectrum of $L=4\frac{n-1}{n-2}\Delta+s$. The third line is obtained from
\begin{align*}
|\phi_i|^\alpha d^*d|\phi_i|^\alpha=\frac{\alpha}{2}|\phi_i|^{2\alpha-2}d^*d|\phi_i|^2-\alpha(\alpha-2)|\phi_i|^{2\alpha-2}|d|\phi_i||^2.
\end{align*}
Using \[\frac{1}{2}d^*d\langle \phi_i,\phi_i\rangle=\langle\nabla^*\nabla\phi_i,\phi_i\rangle-|\nabla\phi_i|^2=\langle D^2\phi_i,\phi_i\rangle-\frac{s}{4}|\phi_i|^2-|\nabla\phi_i|^2\]
and 
\[|\nabla^\lambda\phi_i|^2=|\nabla\phi_i|^2-2\Re\frac{\lambda}{n}\langle (D-\lambda)\phi_i,\phi_i\rangle-\frac{\lambda^2}{n}|\phi_i|^2,\]
we get
\begin{align*}
\left(\frac{\mu}{4}-\frac{n-1}{n}\lambda^2\right)& \Vert |\phi_i|^{\alpha}\Vert^2 \leq\int\limits |\phi_i|^{2\frac{n-2}{n-1}-2}\left(\frac{n}{n-1} |d|\phi_i||^2-|\nabla^\lambda\phi_i|^2\right)\vo_g\\
&+\int\limits |\phi_i|^{2\frac{n-2}{n-1}-2}\langle (D^2-\lambda^2)\phi_i,\phi_i\rangle\vo_g\\
&- \int\limits 2 |\phi_i|^{2\frac{n-2}{n-1}-2}\Re\frac{\lambda}{n}\langle (D-\lambda)\phi_i,\phi_i\rangle \vo_g\\
&\leq\int\limits |\phi_i|^{2\frac{n-2}{n-1}-2}\left(\frac{n}{n-1} |d|\phi_i||^2-|\nabla^\lambda\phi_i|^2\right)\vo_g\\
&\phantom{==}+\int\limits |\phi_i|^{2\frac{n-2}{n-1}-2}\langle(D-\lambda)^2\phi_i,\phi_i\rangle \vo_g\\
&\phantom{==}+ \int\limits 2\left(1-\frac{1}{n}\right)\lambda |\phi_i|^{2\frac{n-2}{n-1}-2}\Re\langle (D-\lambda)\phi_i,\phi_i\rangle \vo_g.
\end{align*}
The last two summands vanish in the limit $i\to \infty$ since $\Vert \phi_i\Vert_{2\frac{n-3}{n-1}}\leq \Vert \phi_i\Vert_2=1$ and, thus, 
\[\left|\int\limits |\phi_i|^{2\frac{n-2}{n-1}-2}\langle (D-\lambda)^2\phi_i,\phi_i\rangle \vo_g\right|\leq \Vert (D-\lambda)^2\phi_i\Vert\ \Vert\, |\phi_i|^{\frac{n-3}{n-1}}\Vert\to 0\]
and
\[ \Big|\int\limits |\phi_i|^{2\frac{n-2}{n-1}-2}\Re\langle (D-\lambda)\phi_i,\phi_i\rangle\vo_g\Big|\leq \Vert (D-\lambda)\phi_i\Vert\ \Vert\, |\phi_i|^{\frac{n-3}{n-1}}\Vert\to 0.\]

For the other -- the first -- summand we use the Kato-type inequality of Lemma \ref{reKI2}
\[ |d|\psi||\leq |(D-\lambda)\psi|+k|\nabla^\lambda\psi|\] which holds outside the zero set of $\psi$. Due to \cite[(3.9)]{CGH} we have $k=\sqrt{\frac{n-1}{n}}$. Thus, we can estimate 
\begin{align*} \int\limits |\phi_i&|^{2\frac{n-2}{n-1}-2}\left(\frac{n}{n-1} |d|\phi_i||^2-|\nabla^\lambda\phi_i|^2\right)\vo_g\\
&= \int\limits |\phi_i|^{2\frac{n-2}{n-1}-2}(k^{-1} |d|\phi_i||-|\nabla^\lambda\phi_i|)(k^{-1} |d|\phi_i||+|\nabla^\lambda\phi_i|)\vo_g\\
&\leq k^{-1}\hspace{-0.2cm}\int\limits_{\{|d|\phi_i|| \geq k|\nabla^\lambda\phi_i|\}} |\phi_i|^{2\frac{n-2}{n-1}-2}|(D-\lambda)\phi_i|(k^{-1} |d|\phi_i||+|\nabla^\lambda\phi_i|)\vo_g\\
&\leq 2k^{-2}\hspace{-0.2cm}\int\limits_{\{|d|\phi_i||\geq k|\nabla^\lambda\phi_i|\}} |\phi_i|^{2\frac{n-2}{n-1}-2}|(D-\lambda)\phi_i| |d|\phi_i||\vo_g\\
&\leq 2k^{-2} \int\limits \left({2\frac{n-2}{n-1}-1}\right)^{-1}|(D-\lambda)\phi_i| |d|\phi_i|^{2\frac{n-2}{n-1}-1}|\vo_g\\
&\leq 2k^{-2}\frac{n-1}{n-3} \Vert (D-\lambda)\phi_i\Vert\ \Vert d|\phi_i|^{\frac{n-3}{n-1}}\Vert.
\end{align*}

For $n\geq 5$ we have $1\geq\frac{n-3}{n-1}\geq \frac{1}{2}$ and, thus, $\Vert d|\phi_i|^{\frac{n-3}{n-1}}\Vert$ is bounded. Together with $\Vert (D-\lambda)\phi_i\Vert\to 0$ we obtain the following: For all $\epsilon>0$ there is an $i_0$ such that for all $i\geq i_0$ we have
\[\int\limits |\phi_i|^{2\frac{n-2}{n-1}-2}\left(\frac{n}{n-1} |d|\phi_i||^2-|\nabla^\lambda\phi_i|^2\right)\vo_g\leq \epsilon.\]
Thus, for $i\to \infty$ the right handside of inequality \eqref{ineq-ess} tends to zero. It remains to show that on the left handside of this inequality $\Vert \phi_i\Vert_{2\alpha}$ is bounded from below:

Assume that $\Vert \phi_i\Vert_{2\alpha}\to 0$ as $i\to \infty$. Then from \eqref{non-zero-in} with $\beta=\frac{n-2}{n-1}=\alpha$ and the scalar curvature $s$ bounded from below, we see that then $\Vert d |\phi_i|^\alpha\Vert\to 0$ as $i\to \infty$ and moreover, \eqref{non-zero-in} then reads 

\begin{align*}
0&\leftarrow \Vert (D-\lambda)^2\phi_i\Vert_2\\
&\geq (2\alpha-1)\frac{1}{\alpha^2}\int\limits_{|\phi_i|\ne 0} |d|\phi_i|^\alpha|^2\vo_g+ \int\limits\left(\frac{s}{4}-\frac{n-1}{n}\lambda^2\right) |\phi_i|^{2\alpha}\vo_g\\
&\phantom{=}-2\frac{n-1}{n}\lambda \Vert(D-\lambda)\phi_i\Vert\\
&= \frac{1}{4}\left( 4\frac{n-1}{n-2}\int |d|\phi_i|^\alpha|^2\vo_g+ \int s |\phi_i|^{2\alpha}\vo_g\right)
-\frac{n-1}{n}\lambda^2 |\phi_i|^{2\alpha}\vo_g\\
 &\phantom{=}-\frac{n-1}{(n-2)^2}\int |d|\phi_i|^\alpha|^2\vo_g -2\frac{n-1}{n}\lambda \Vert(D-\lambda)\phi_i\Vert\\
&= \frac{1}{4}\int |\phi_i|^\alpha L_g |\phi_i|^\alpha|\vo_g -\frac{n-1}{(n-2)^2}\int |d|\phi_i|^\alpha|^2\vo_g \\
&\phantom{=}- \frac{n-1}{n}\lambda^2 |\phi_i|^{2\alpha}\vo_g-2\frac{n-1}{n}\lambda \Vert(D-\lambda)\phi_i\Vert
\end{align*} 

which implies that $\int |\phi_i|^\alpha L_g |\phi_i|^\alpha|\vo_g\to 0$ as $i\to \infty$ and contradicts $Q>0$.
Hence, $\Vert \phi_i\Vert_{2\alpha}$ is bounded from below and in the limit for $i\to\infty$ inequality \eqref{ineq-ess} leads to $\frac{\mu}{4}\leq \frac{n-1}{n}\lambda^2$.
\end{proof}

\begin{remark}
We think that the restriction on the dimension and the assumption that the scalar curvature is bounded from below could be circumvented and just appear in our proof for technical reasons. But we cannot get rid of these assumptions until now. Especially the scalar curvature assumption one would expect to be not necessary since in the extremal case, when the scalar curvature goes to $-\infty$ in every direction, we already have $\mu<0$. But in general there could just be a tiny neighbourhood of a ray where $s\to -\infty$ and our proof cannot handle this up to now.  
\end{remark}

With the two theorems on Hijazi inequalities from above, we can now prove the conformal Hijazi inequality:

\begin{proof}[Theorem \ref{Hij2}] 
For $Q\leq 0$ the inequality is trivially satisfied. Thus, we restrict ourselves to the case $Q> 0$:\\
We may assume that $g$ is itself a complete metric of finite volume satisfying the condition (i): $0\not\in\sigma_{ess}(D_g)$. Due to Lemma \ref{sequ} there exists a sequence $g_i$ of complete metrics of unit volume with $g_i\equiv g$ near infinity and $\lambda_1^+(g_i)\to \lm$.\\
We first consider the case that there is an infinite subsequence $g_{i_j}$ such that $\lambda_1^+(g_{i_j})$ is an eigenvalue of $D_{g_{i_j}}$. Then we can apply Theorem \ref{Hij1} and equality \eqref{Qinf} and obtain \[\lambda_1^+(M,g_{i_j},\sigma)^2\geq \frac{n}{4(n-1)} \mu (M,g_{i_j}) \geq \frac{n}{4(n-1)} Q(M,g).\]
Thus, for $j\to \infty$ we obtain the conformal Hijazi inequality.\\
Now we consider the remaining case -- only finitely many $\lambda_1^+(g_i)$ are eigenvalues. Thus, from Lemma \ref{spec} we know that then there is an infinite subsequence $g_{i_j}$ such that $\lambda_1^+(g_{i_j})\in \sigma_{ess}(D_{g_{i_j}})$. But if for two metrics $g_i$ and $g_k$ we have $\sigma_{ess}(D_{g_i})\ni\lambda_1^+(g_i)\geq \lambda_1^+(g_k)\in \sigma_{ess}(D_{g_k})$, then  $\lambda_1^+(g_i)$ already equals $\lambda_1^+(g_k)$ since $g_k\equiv g_i$ near infinity and the essential spectrum only depends on the manifold at infinity. Hence, there has to exist a constant subsequence $\lm=\lambda_1^+(g_{i_j})\in \sigma_{ess}(D_{g_{i_j}})=\sigma_{ess}(D_g)$. Lemma \ref{zero}.iii then gives $\lm=0$ and, thus, $0\in\sigma_{ess}(D_g)$. This is a contradiction to the assumption.\par
So we assume now that $0\in \sigma_{ess}(D)$. Then condition (ii) has to be fulfilled and Theorem \ref{n5sb} implies $\mu\leq 0$ and, thus, $Q\leq 0$.
\end{proof}

\begin{example}
We consider the Riemannian manifold $(M \times \R, g_M+dt^2)$ where $(M, g_M)$ is closed, spin and has positive scalar curvature. Due Example \ref{rnconf} $M\times \R$ is conformally parabolic and the conformal metric $\overline{g}=f(t)^2g$ where $f$ is positive and smooth and  $f(t)=\frac{1}{t^2}$ for $|t|\geq 1$ is complete and of finite volume. Its scalar curvature is bounded from below. Note that $Q(M\times \R, g_m+dt^2)>0$ for $M$ having positive scalar curvature \cite[Prop. 5.7]{Aku}.

Then with Theorem \ref{Hij2} we know that at least for $n\geq 5$ the conformal Hijazi inequality is valid. Furthermore, for $n=3$ and $n=4$  and $0\not\in \sigma_{ess}D^{\overline{g}}$ Theorem \ref{Hij2} also gives the validity of the conformal Hijazi inequality. In case that $0\not\in \sigma_{ess}(D^{\overline{g}})$ we know from $\lm=0$. Thus, if such a two- or three-dimensional manifold $M$ exists such that $0\in \sigma_{ess}(D^{\overline{g}})$ the conformal inequality does not hold for this manifold. But, we don't know yet whether such a manifold exists.
\end{example}

\textit{Acknowledgement} I am grateful to Hans-Bert Rademacher for discussions, help and his patient support. I also thank the MPI for Mathematics in the Sciences for financial support.\\

\bibliography{Hij}
\bibliographystyle{plain}

\textsc{Universit\"at Leipzig, Fakult\"at f\"ur Mathematik und Informatik, Johannisgasse 26, Leipzig 04013, Germany;\\ nadine.grosse@mathematik.uni-leipzig.de}  

\end{document}